\documentclass[amssymb,amstex,amsmath,10pt,a4paper]{amsart}
\usepackage{amsmath,amssymb,latexsym,amsfonts,url}

\newtheorem{thm}{Theorem}
\newtheorem*{theorem}{Theorem}
\newtheorem{df}[thm]{Definition}

\newtheorem{cor}[thm]{Corollary}
\newtheorem{lem}[thm]{Lemma}

\newcommand{\vol}{{\rm Vol}}

\newcommand{\Ric}{{\rm Ric}}

\newenvironment{pf}{\noindent {\it Proof.}}{\\}

\numberwithin{equation}{section}

\begin{document}
\title[$L^2$ harmonic 1-forms on minimal submanifolds]{$L^2$ harmonic 1-forms on minimal submanifolds \\in hyperbolic space}
\author{Keomkyo Seo}
\address{Department of Mathematics\\
Sookmyung Women's University\\
Hyochangwongil 52, Yongsan-ku \\
Seoul, 140-742, Korea}
\email{kseo@sookmyung.ac.kr}

\thanks{This research was supported by the Sookmyung Women's University Research Grants 2009.}

\begin{abstract}
In this paper, we prove the nonexistence of $L^2$ harmonic 1-forms on a complete super stable minimal submanifold $M$ in hyperbolic space under the assumption that the first eigenvalue $\lambda_1 (M)$ for the Laplace operator on $M$ is bounded below by $(2n-1)(n-1)$. Moreover, we provide sufficient conditions for minimal submanifolds in hyperbolic space to be super stable. \\

\end{abstract}

\maketitle

%%%%%%%%%%%%%%%%%%%%%%%%%%%%%%%%%%%%%%%%%%%%%%%%
\section{Introduction}
%%%%%%%%%%%%%%%%%%%%%%%%%%%%%%%%%%%%%%%%%%%%%%%%
Let $M$ be an $n$-dimensional complete minimal submanifold in $\mathbb{R}^{n+m}$. In \cite{Miyaoka}, Miyaoka showed that if $M$ is a complete stable minimal hypersurface in $\mathbb{R}^{n+1}$, then there are no nontrivial $L^2$ harmonic $1$-forms on $M$. Recall that a minimal submanifold is said to be $stable$ if the second variation of its volume is always nonnegative for any normal variation with compact support.  Yun \cite{Yun} proved that if  $M$ is a complete minimal hypersurface with $\displaystyle{\Big(\int_M |A|^n dv \Big)^{1 \over n} < C_2 = \sqrt{C_s^{-1}}}$, then there are no nontrivial $L^2$ harmonic $1$-forms on $M$. Here $C_s$ is a Sobolev constant in \cite{MS}. Recently, the author extended this result to higher codimension cases in a Euclidean space and hyperbolic space \cite{Seo08-1, Seo10}. We denote by $-|A|^2$ and $\int_M |A|^n dv$ the {\it scalar curvature} and  the {\it total scalar curvature} of a minimal submanifold $M$ in a Euclidean space, respectively.

This paper is concerned with the nonexistence of $L^2$ harmonic 1-forms on complete minimal submanifolds in hyperbolic space. We shall denote by
$\mathbb{H}^n$ the $n$-dimensional hyperbolic space of constant sectional curvature $-1$. For an $n$-dimensional minimal submanifold $M$ in $\mathbb{H}^{n+m}$, a simple computation shows that for a variational vector field $E=\varphi \nu$, the second variation of its volume $\vol (M_t)$ satisfies \cite{Lawson}
\begin{eqnarray*}
\frac{d^2 \vol(M_t)}{dt^2} \geq \int_M |\nabla \varphi|^2 - (|A|^2-n)\varphi^2 dv,
\end{eqnarray*}
where $\varphi \in W_0^{1,2}(M)$ and $\nu$ is the unit normal vector field. (See also \cite{Spruck}.) Motivated by this, we introduce the concept of super stability in hyperbolic space as follows:
\begin{df}
\rm{We call a minimal submanifold $M$ in $\mathbb{H}^{n+m}$ $super
\ stable$ if for any $\varphi \in W_0^{1,2}(M)$}
\begin{eqnarray*}
\int_M |\nabla \varphi|^2  - (|A|^2-n) \varphi^2 dv \geq 0 .
\end{eqnarray*}
\end{df}
Note that when $m=1$, the concept of super stability is the same as the usual definition of stability. It should be mentioned that Q. Wang \cite{Wang} introduced the concept of super stability for minimal submanifolds in a Euclidean space and proved that the only complete super stable minimal submanifold $M^n (n \geq 3)$ of finite total scalar curvature in a Euclidean space is an affine $n$-plane. Recently the author \cite{Seo08-2} proved that if $M$ is a complete immersed super stable minimal submanifold in a Euclidean space with flat normal bundle satisfying that $\int_M |A|^2 < \infty$, then $M$ is an affine $n$-plane.

The first eigenvalue for the Laplace operator $\Delta$ on a Riemannian manifold $M$ is defined as
\begin{align*}
\lambda_1 (M) = \inf_{f} \frac{\int_M |\nabla f|^2}{\int_M f^2} ,
\end{align*}
where the infimum is taken over all compactly supported smooth functions on $M$. It was proved by Cheung and Leung
\cite{Cheung and Leung}  that for an $n$-dimensional complete minimal submanifold $M$ in $\mathbb{H}^{n+m}$, the first eigenvalue $\lambda_1 (M)$ satisfies
\begin{eqnarray} \label{ineq:CL}
\frac{1}{4}(n-1)^2 \leq \lambda_1 (M) .
\end{eqnarray}
Here this inequality is sharp because equality holds when $M$ is
totally geodesic \cite{McKean}. In Section 2, we shall prove a gap theorem for $L^2$ harmonic 1-forms on a complete super stable minimal submanifold in hyperbolic space. More precisely, we prove
\begin{theorem}
Let $M$ be a complete super stable minimal submanifold in $\mathbb{H}^{n+m}$. Assume that the first eigenvalue of $M$ satisfies
\begin{eqnarray*}
(2n-1)(n-1) < \lambda_1 (M) .
\end{eqnarray*}
Then there are no nontrivial $L^2$ harmonic 1-forms on $M$.
\end{theorem}

For a harmonic function $f$ on a complete Riemannian manifold with finite Dirichlet energy, it is well known that its differential $df$ is $L^2$ harmonic $1$-form on $M$. The topology of a stable minimal hypersurface in a nonnegatively curved manifold is closely related with the space of bounded harmonic functions with finite Dirichlet energy.(\cite{CSZ}, \cite{LW}, \cite{Wei}) In this direction, our main theorem provides the nonexistence of nontrivial harmonic functions on a complete super stable minimal submanifold in hyperbolic space with finite Dirichlet energy under the first eigenvalue assumption.

In Section 3, we find sufficient conditions for minimal submanifolds in hyperbolic space to be super stable.

%%%%%%%%%%%%%%%%%%%%%%%%%%%%%%%%%%%%%%%%%%%%%%%%%%%%%%%%%%%%%%%%%%%%%%%%
\section{A gap theorem for $L^2$ harmonic 1-forms}
%%%%%%%%%%%%%%%%%%%%%%%%%%%%%%%%%%%%%%%%%%%%%%%%%%%%%%%%%%%%%%%%%%%%%%%%
We begin with the following lemma.
\begin{lem} {\rm (\cite{Leung})}  \label{lem:Leung}
Let $M$ be an $n$-dimensional complete immersed minimal submanifold in
$\mathbb{H}^{n+m}$. Then the Ricci curvature of $M$ satisfies
\begin{eqnarray*}
\Ric(M) \geq -(n-1) -{n-1 \over n}|A|^2.
\end{eqnarray*}
\end{lem}
We are now ready to prove our main theorem.
\begin{thm} \label{thm:main}
Let $M$ be a complete super stable minimal submanifold in $\mathbb{H}^{n+m}$. Assume that the first eigenvalue of $M$ satisfies
\begin{eqnarray*}
(2n-1)(n-1) < \lambda_1 (M).
\end{eqnarray*}
Then there are no nontrivial $L^2$ harmonic 1-forms on $M$.
\end{thm}
\begin{pf}
Let $\omega$ be an $L^2$ harmonic 1-form on M, i.e.,
\begin{eqnarray*}
\Delta \omega = 0 \ \ {\rm and}\ \  \int_M |\omega|^2 dv < \infty .
\end{eqnarray*}
In an abuse of notation, we will refer to a harmonic $1$-form and its dual harmonic vector field both by $\omega$. Bochner
formula says
\begin{eqnarray*}
\Delta|\omega|^2 = 2(|\nabla \omega|^2 + \Ric(\omega,\omega)) .
\end{eqnarray*}
Moreover, we have
\begin{eqnarray*}
\Delta|\omega|^2 = 2(|\omega|\Delta |\omega| + |\nabla|\omega||^2) .
\end{eqnarray*}
Thus using Kato type inequality\cite{X.Wang} we get
\begin{eqnarray*}
|\omega| \Delta |\omega| - \Ric(\omega,\omega) = |\nabla \omega|^2
-|\nabla|\omega||^2 \geq {1 \over {n-1}} |\nabla|\omega||^2 .
\end{eqnarray*}
Therefore applying Lemma \ref{lem:Leung} gives
\begin{eqnarray} \label{ineq:1-form}
|\omega| \Delta |\omega|  + {{n-1} \over n}|A|^2 |\omega|^2 + (n-1) |\omega|^2 \geq {1 \over {n-1}} |\nabla|\omega||^2 .
\end{eqnarray}
Furthermore, the super stability of $M$ implies
\begin{eqnarray*}
\int_M |\nabla \phi|^2 - (|A|^2 -n)\phi^2 dv \geq 0
\end{eqnarray*}
all $\phi \in W_0^{1,2}(M)$ .
Substituting $|\omega|\phi$ for $\phi$ gives
\begin{eqnarray*}
\int_{M} |\nabla (|\omega|\phi)|^2 - (|A|^2 -n)|\omega|^2\phi^2 dv\geq 0.
\end{eqnarray*}
Thus applying integration by parts and divergence theorem we have
\begin{align} \label{ineq:test}
0 &\leq -\int_M |\omega|\phi\Delta (|\omega|\phi) dv - \int_M (|A|^2 -n)|\omega|^2\phi^2 dv \nonumber \\
  &= -\int_M |\omega|\phi (|\omega|\Delta \phi + \phi \Delta |\omega| + 2 \langle\nabla|\omega|, \nabla \phi \rangle) dv - \int_M (|A|^2 -n)|\omega|^2\phi^2 dv \nonumber \\
  &= -\int_M \phi |\omega|^2 \Delta \phi dv - \int_M \phi^2 (|\omega|\Delta|\omega| +  |A|^2 |\omega|^2) dv \nonumber \\
   &\quad - 2\int_M |\omega|\phi \langle\nabla|\omega|, \nabla \phi \rangle dv + n \int_M |\omega|^2 \phi^2 dv \nonumber \\
  &=\int_M \langle \nabla (\phi|\omega|^2), \nabla \phi \rangle dv  - \int_M \phi^2 (|\omega|\Delta|\omega| +  |A|^2 |\omega|^2) dv \nonumber \\
  &\quad - 2\int_M |\omega|\phi \langle\nabla|\omega|, \nabla \phi \rangle dv + n \int_M |\omega|^2 \phi^2 dv \nonumber \\
  &=\int_M |\omega|^2|\nabla \phi|^2 dv - \int_M \phi^2 (|\omega|\Delta|\omega| +  |A|^2 |\omega|^2) dv + n \int_M |\omega|^2 \phi^2 dv .
\end{align}
Combining the inequalities (\ref{ineq:1-form}) and (\ref{ineq:test}) gives
\begin{align} \label{ineq:last}
0 &\leq (2n-1)\int_M |\omega|^2 \phi^2 dv + \int_M |\omega|^2|\nabla \phi|^2 dv -\frac{1}{n-1}\int_M |\nabla|\omega||^2 \phi^2 dv \nonumber \\
  &\quad - \frac{1}{n}\int_M |A|^2 |\omega|^2 \phi^2 dv.
\end{align}
On the other hand, since the first eigenvalue satisfies
\begin{align*}
\lambda_1 (M) \leq \frac{\int_M |\nabla \phi|^2 dv}{\int_M |\phi|^2 dv}
\end{align*}
for all $\phi \in W_0^{1,2}(M)$, we substitute $|\omega|\phi$ for $\phi$. Then we get
\begin{align} \label{ineq:lambda_1}
\lambda_1 (M)\int_M |\omega|^2 \phi^2 dv \leq \int_M |\omega|^2|\nabla \phi|^2 + |\phi|^2|\nabla |\omega||^2  + 2\phi|\omega| \langle \nabla \phi, \nabla |\omega| \rangle dv.
\end{align}
It follows from the inequalities (\ref{ineq:last}) and (\ref{ineq:lambda_1}) that
\begin{align*}
0 &\leq \Big(\frac{2n-1}{\lambda_1 (M)} + 1\Big) \int_M |\omega|^2|\nabla \phi|^2 dv + \Big(\frac{2n-1}{\lambda_1 (M)} -\frac{1}{n-1} \Big)\int_M |\nabla|\omega||^2 \phi^2 dv \\
  &\quad + \frac{2(2n-1)}{\lambda_1 (M)} \int_M \phi|\omega| \langle \nabla \phi, \nabla |\omega| \rangle dv - \frac{1}{n}\int_M |A|^2 |\omega|^2 \phi^2 dv.
\end{align*}
Combining the Cauchy-Schwarz inequality with Young's inequality for $\varepsilon >0$ yield
\begin{align*}
 2\int_M \phi|\omega| \langle \nabla \phi, \nabla |\omega| \rangle dv \leq \frac{\varepsilon}{2}\int_M |\nabla|\omega||^2\phi^2 dv + \frac{2}{\varepsilon} \int_M |\omega|^2|\nabla \phi|^2 dv ,
\end{align*}
Thus we have
\begin{align*}
0 &\leq \Big\{\Big(1+\frac{2}{\varepsilon}\Big)\frac{2n-1}{\lambda_1 (M)} + 1\Big\} \int_M |\omega|^2|\nabla \phi|^2 dv \\
  &\quad + \Big\{\Big(1+\frac{\varepsilon}{2}\Big)\frac{2n-1}{\lambda_1 (M)} -\frac{1}{n-1} \Big\}\int_M |\nabla|\omega||^2 \phi^2 dv - \frac{1}{n}\int_M |A|^2 |\omega|^2 \phi^2 dv.
\end{align*}
Now fix a point $p \in M$ and for $R>0$ choose a cut-off function
satisfying $0 \leq \phi \leq 1$, $\phi \equiv 1$ on
$B_p(R)$, $\phi = 0 $ on $M \setminus B_p(2R)$, and
$\displaystyle{|\nabla \phi| \leq {1 \over R}}$. Since $\lambda_1 (M) > (2n-1)(n-1)$ by assumption, choosing $\varepsilon
> 0$ sufficiently small and letting $R \rightarrow \infty$, we
obtain $|\omega| \equiv 0$ or $|A| \equiv 0$. If $|\omega| \equiv 0$, then there is no nontrivial $L^2$ harmonic 1-form on $M$. If $|A| \equiv 0$, $M$ is isometric to $\mathbb{H}^n$. However, it was proved by Dodziuk \cite{Dodziuk1} that there is no nontrivial $L^2$ harmonic 1-form on $\mathbb{H}^n$. Hence we get the conclusion. \qed
\end{pf}

As mentioned in the introduction, when $m=1$, the definition of super stability is the same as that of stability. Thus it immediately follows
\begin{cor}
Let $M$ be a complete stable minimal hypersurface in $\mathbb{H}^{n+1}$ satisfying that
\begin{eqnarray*}
(2n-1)(n-1) < \lambda_1 (M) .
\end{eqnarray*}
Then there are no nontrivial $L^2$ harmonic 1-forms on $M$.
\end{cor}
It was shown by Palmer \cite{Palmer} that if there exists a codimension one cycle $C$ in a complete minimal hypersurface $M$ in $\mathbb{R}^{n+1}$ which  does not separate $M$, then $M$ is unstable. We prove an analogue of his result in hyperbolic space.
\begin{cor}
Let $M$ be a complete minimal submanifold in $\mathbb{H}^{n+m}$ satisfying that
\begin{eqnarray*}
(2n-1)(n-1) < \lambda_1 (M).
\end{eqnarray*}
Suppose there exists a codimension one cycle $C$ in $M$ which does not separate $M$. Then $M$ cannot be super stable.
\end{cor}
\begin{pf}
By the results of Dodziuk \cite{Dodziuk2}, there exists a nontrivial $L^2$ harmonic 1-form on $M$. However this is impossible by Theorem \ref{thm:main}. \qed
\end{pf}

Let $f$ be a harmonic function on $M$ with finite Dirichlet energy, i.e., $\Delta f = 0$ and $\int_M |\nabla f|^2 dv < \infty$. Then its differential $df$ is $L^2$ harmonic 1-form on $M$. Hence we get
\begin{cor}
Let $M$ be a complete super stable minimal submanifold in $\mathbb{H}^{n+m}$ satisfying that
\begin{eqnarray*}
(2n-1)(n-1) < \lambda_1 (M).
\end{eqnarray*}
Then there are no nontrivial harmonic functions on $M$ with finite Dirichlet energy.
\end{cor}
%%%%%%%%%%%%%%%%%%%%%%%%%%%%%%%%%%%%%%%%%%%%%%%%%%%%%%%%%%%%%%%%%%%%%%%%
\section{Super stability of minimal submanifolds in hyperbolic space}
%%%%%%%%%%%%%%%%%%%%%%%%%%%%%%%%%%%%%%%%%%%%%%%%%%%%%%%%%%%%%%%%%%%%%%%%
In this section, we provide two sufficient conditions for minimal submanifolds in hyperbolic space to be super stable.
\begin{thm}
Let $M$ be a complete minimal submanifold in $\mathbb{H}^{n+m}$.
If $|A| \leq \frac{(n+1)^2}{4}$ at every point in $M$, then $M$ is
super stable.
\end{thm}
\begin{pf}
Using the inequality (\ref{ineq:CL}), we have
\begin{eqnarray*}
\frac{(n-1)^2}{4} \leq \lambda_1 (M) \leq \frac{\int_M |\nabla
f|^2 dv}{\int_M f^2 dv}
\end{eqnarray*}
for all $f\in W_0^{1,2}(M)$. Hence
the assumption that $|A|^2 \leq \frac{(n+1)^2}{4}$ implies
\begin{eqnarray*}
\int_M |\nabla f|^2 - (|A|^2 - n)f^2 dv \geq \int_M (\lambda_1 (M) +
n - |A|^2 )f^2 dv \geq 0,
\end{eqnarray*}
which completes the proof. \qed
\end{pf}

It is well known that the following Sobolev inequality \cite{HS} on an $n$-dimensional minimal submanifold $M$ in $\mathbb{H}^{n+m} (n\geq 3)$ holds
\begin{eqnarray} \label{ineq:Sobolev}
\Big(\int_M |f|^{\frac{2n}{n-2}} dv\Big)^{\frac{n-2}{n}} \leq C_s
\int_M |\nabla f|^2 dv  \text{\hspace{0.5cm} for } f \in W_0^{1,2}(M),
\end{eqnarray}
where $C_s$ is the Sobolev constant in \cite{HS}. We shall make use of this inequality to prove the following theorem, which gives another sufficient condition for super stability of minimal submanifolds in hyperbolic space.
\begin{thm}
Let $M$ be an $n$-dimensional complete minimal submanifold in $\mathbb{H}^{n+m}, n \geq 3$.
If $\int_M |A|^n dv \leq (\frac{1}{C_s})^{\frac{n}{2}}$, then $M$
is super stable.
\end{thm}
\begin{pf}
It is sufficient to show that
\begin{eqnarray*}
\int_M |\nabla f|^2 - (|A|^2 - n)f^2 dv \geq 0
\end{eqnarray*}
for all $f\in W_0^{1,2}(M)$. By Sobolev
inequality (\ref{ineq:Sobolev}), we have
\begin{equation} \label{ineq:1}
\int_M |\nabla f|^2 - (|A|^2 - n)f^2 dv \geq \frac{1}{C_s}
\Big(\int_M |f|^{\frac{2n}{n-2}} dv\Big)^{\frac{n-2}{n}} - \int_M
|A|^2 f^2 dv .
\end{equation}
On the other hand, it follows from H\"{o}lder inequality that
\begin{eqnarray} \label{ineq:2}
\int_M |A|^2 f^2 dv \leq \Big(\int_M |A|^n dv \Big)^{\frac{2}{n}}
\Big(\int_M |f|^{\frac{2n}{n-2}} dv \Big)^{\frac{n-2}{n}}.
\end{eqnarray}
Combining (\ref{ineq:1}) with (\ref{ineq:2}) we have
\begin{eqnarray*}
\int_M |\nabla f|^2  - (|A|^2 - n)f^2 dv &\geq&
\Big\{\frac{1}{C_s} - \Big(\int_M |A|^n
dv\Big)^{\frac{2}{n}}\Big\} \Big(\int_M |f|^{\frac{2n}{n-2}} dv
\Big)^{\frac{n-2}{n}} \\
 &\geq& 0,   \hspace{1cm} \ \mbox{(by assumption)}
\end{eqnarray*}
which completes the proof. \qed
\end{pf}

\subsection*{Acknowledgment}The author would like to thank the referee for the helpful comments and suggestions.

\end{document}